\documentclass[a4paper, abstracton]{scrarticle}

\usepackage[utf8]{inputenc}
\usepackage{lmodern}
\usepackage[T1]{fontenc}
\usepackage[english]{babel}
\usepackage{amsmath,amssymb, amsthm}

\usepackage{graphicx}
\usepackage{url}
\usepackage{booktabs}
\usepackage{hhline}
\usepackage{subcaption}

\usepackage{authblk}

\begin{document}

\title{Sparse Linear Surrogates for Interpretable Budget Allocation}

\author[1]{Marc Goerigk}
\author[1]{Michael Hartisch}
\author[1]{Sebastian Merten}

\affil[1]{Business Decisions and Data Science, University of Passau,\authorcr Dr.-Hans-Kapfinger-Str. 30, 94032 Passau, Germany}

\date{}

\maketitle

\begin{abstract}
To address the demand for inherently interpretable optimization methods, we introduce novel linear surrogates for budget allocation problems. These surrogates consist of sparse linear rules that map instances to feature-based representations of solutions. We present an exact approach based on mixed-integer programming as well as a heuristic for their computation. The performance of both approaches is analyzed through computational experiments.
\end{abstract}

\section{Introduction}
\label{sec:intro}

The growing use of automated decision-support tools has increased the need for transparency and comprehensibility. The ability to explain is particularly important when optimized solutions affect workers, planners, customers, or other stakeholders not involved in the decision process. When these stakeholders cannot understand the logic behind a decision, trust and acceptance may decline, limiting the practical value of optimized solutions. Hence, providing easily comprehensible insights into when and why a solution is recommended should be considered as an important aspect when designing a decision-support tool.
We argue that interpretability---where not just a post-hoc justification of an obtained result is provided, but the solution process is made accessible---should receive greater attention in this context. Utilizing an existing framework for finding interpretable surrogates of the optimization process \cite{GOERIGK20231312}, we investigate sparse linear surrogates. In this data-based approach, we aim at finding a linear rule that maps instances--represented by feature vectors---to feature-based representations of solutions \cite{goerigk2024feature}, i.e., a vector of instance features is mapped to a vector of solution features. This differs from previous approaches, where only decision trees where considered as the implemented surrogate,
and is in line with the observation that no one-size-fits-all explanation can exist, making the availability of various approaches indispensible \cite{aydin2026explainable}. Such a linear surrogate can then be used both as a surrogate for an optimization process and as a transparent policy guideline. In this paper, we study such an approach in the setting of the budget allocation problem, which we formulate as a knapsack problem.

Our work relates to interpretability and explainability in optimization \cite{aydin2026explainable},
and contextual optimization \cite{elmachtoub2022smart}. For a broader overview, we refer to \cite{goerigk2024feature}. Our use of linear surrogates is also related to linear decision rules in robust optimization \cite{kuhn2011primal}, for which sparse optimal representations exist in relevant problem classes \cite{lu2026sparsity}. There, however, decision rules are primarily used to obtain computationally tractable approximations of adaptive policies. Finally, feature-based approaches, particularly for the newsvendor problem, use regularized regression models to map instance features directly to decisions \cite{HUBER2019904,ban2019big}.
In contrast, our method learns sparse linear surrogates that map instance features to so-called meta-solutions: feature-based descriptions of restricted solution spaces from which implementable solutions can subsequently be recovered. 

As an example, consider a monthly budget-allocation problem in which a factory divides its budget between machine-maintenance tasks and product-innovation projects. In this simple example, we consider a single instance feature $f\in[0,1]$, representing the share of production activities disrupted by machine failures. An implementable solution selects individual maintenance tasks and innovation projects, whereas the meta-solution records only the fraction $b\in[0,1]$ of the budget allocated to maintenance. Thus, $b$ does not determine the actual selection but provides an interpretable restriction of the solution space.
Suppose historical data indicate that the maintenance share associated with the best outcomes is approximately $b=0.6f+0.2$. Hence, the maintenance share should increase gradually as the disruption rate increases.
Obviously, a linear surrogate represents this relationship exactly. In contrast, a decision tree can return only as many meta-solutions as there are leaves. For example, a depth-one tree with a split at $f\leq 0.5$ assigns instances with disruption rates $f=0.1$ and $f=0.5$ to the same meta-solution despite requiring different allocations. A deeper tree provides more levels but still approximates the gradual relationship by a small number of steps.

\section{Finding Sparse Linear Surrogates}
\label{sec:problem}
Following the notation used in \cite{goerigk2024feature} let $F_I \in \mathbb{N}$ and $F_S \in \mathbb{N}$ be the number of instance and solution features, respectively. The instance feature space is given by $\Phi_I \subseteq \mathbb{R}^{F_I}$ and the solution feature space is given by $\Phi_S \subseteq \mathbb{R}^{F_S}$. We call $\phi_I \in \Phi_I$ an instance feature vector and $\phi_S \in \Phi_S$ a meta-solution. We assume the existence of a data base consisting of $N$ instances, described by $\phi_I^1,\ldots,\phi_I^N$. The set of admissible surrogates is given by
$$\mathcal{A}=\left\lbrace a:\Phi_I\rightarrow \Phi_S \ \middle|\ 
a(\phi_I)=A\phi_I + c \textnormal{, with $c\in \mathbb{R}^{F_S}$,  $A \in M(K)$} \right\rbrace , $$  
where $M(K)$ describes a set of sparse $F_S\times F_I$ matrices and $K\in\mathbb{N}$ is the sparsity parameter. We consider three different types of sparsity:
row sparsity where each row has at most $K$ non-zero elements; column sparsity where at most $K$ columns contain non-zero elements; and global sparsity, where the matrix has at most $K$ non-zero elements in total.
 We want to find a surrogate that performs well on average on the data base. To this end, we assume that, once the surrogate has specified a meta-solution $a(\phi_I^j)$, the best implementable solution $x$ satisfying the features prescribed by $a(\phi_I^j)$ is selected from the instance-specific domain $\mathcal{X}(\phi_I^j)$. Assuming an underlying minimization problem, the optimization problem reads

$$\min_{a \in \mathcal{A}} \sum_{j \in [N]} \min_{\substack{x \in \mathcal{X}(\phi_I^j)\\ \phi_S(x)=a(\phi_I^j)}} c(\phi_I^j,x)\ ,$$
where $c(\phi_I^j,x)$ is the cost arising in instance $j$ if solution $x$ is implemented and $\phi_S(x)$ extracts the solution features of $x$. We use $[N]=\{1,\ldots,N\}$ to denote index sets.

To illustrate this framework, we apply it to the budget allocation problem. We assume that for all instances $j\in [N]$  comprehensible features $\phi_I^j$ are given, including instance parameters profits $p^j$, weights $w^j$, and budget $B^j$, but also contextual or aggregated information such as profit-weight ratios.
We assume the selectable items of each instance to be partitioned into $F_S$ different categories (e.g., 
different spending categories), where $P_s$ is the set of items belonging to category $s\in[F_S]$. Then, a meta-solution specifies the share of the budget used for each category. It is worth noting that the dimension of the underlying problem is not fixed, i.e., the number of items may change from instance to instance. This is enabled by the feature-based description of the meta-solution, that only requires the same features to be used to describe the resulting solution space. For a given sparsity type $M(K)$ our problem can be stated as follows:
\begingroup
\allowdisplaybreaks
\begin{subequations}\label{eq:model}
\begin{align}
    \max \ & \sum_{j \in [N]} (p^{j})^\top x^{j} & \\
    \text{s.t.} \ & A\phi_I^{j} + c = b^{j} & \forall j \in [N] \label{c1} \\
    & \sum_{s \in [F_S]} b_s^{j} \leq 1 & \forall j \in [N] \label{c5}\\
   & \sum_{k \in P_s} w_{k}^{j} x_{k}^{j} \leq b_{s}^{j} {B^j} & \forall j \in [N],  s\in [F_S]\label{c6}\\
       & x_k^{j} \in \{0,1\} & \forall j \in [N], s\in [F_S], k \in P_s \label{c7}\\
    & b^{j} \in [0,1]^{F_S}\ \forall j \in [N],\ A \in M(K),\ c \in \mathbb{R}^{F_S} \span \label{c9}
\end{align}
\end{subequations}
\endgroup
For each instance $j\in[N]$, $b^j$ is the solution feature vector describing the share of the budget used in each category and $x^j$ is the solution implemented adhering to this meta-solution.
Depending on the sparsity structure, $A\in M(K)$ implicitly includes auxiliary binary variables and constraints.

While this formulation ensures every prediction $b^j$ to be in $[0,1]^{F_S}$ for all instances $j\in [N]$, this does not necessarily hold for unseen test data. To ensure the existence of feasible solutions for new data points, the prediction is scaled such that $\sum_{s\in F_S}b^j_s=1$, with potentially negative entries set to zero.

As a heuristic, we propose a two-stage approach, in which the tasks of obtaining good meta solutions and finding a suitable predictor are solved sequentially. First, each individual instance $j \in [N]$ is solved and the corresponding optimal vectors $b^{j}$ are extracted. Secondly, we solve the NP-hard~\cite{natarajan1995sparse} problem of finding a sparse linear regression model
heuristically by using the Lasso method~\cite{tibshirani1996regression} and trim the coefficient matrix to fit our sparsity definitions.

\section{Computational Experiments}
We present computational experiments on synthetic instances to evaluate the performance of the proposed methods. As benchmarks, we compute decision tree surrogates (labeled \texttt{dt}) with a maximum depth of two using the CART-based heuristic presented in \cite{goerigk2024feature}. Additionally, we determine the best single meta-solution for the given training data. This can be viewed as using a linear regression model with $A:=0$ or a decision tree-based surrogate only consisting of a single node.
Solving problem~\eqref{eq:model} is referred to as the exact method and abbreviated as \texttt{ex}, the Lasso-based heuristic is abbreviated as \texttt{heu}.

We consider instances with $16$ items grouped into four categories. We generate data by first sampling three basis scenarios, given by profit-to-weight values from $[0.5,3.5]$ for every item category. For each data point, one basis scenario is randomly selected. For each item we sample a profit uniformly from $[10,30]$. Based on the selected basis scenario, a weight is sampled accordingly from an interval within $\pm 5\%$ of the prescribed profit-to-weight ratio. Finally, all items are scaled such that their weights sum up to 100. The budget is set to 50 for all instances. As instance features we make use of \textit{individual} instance features, consisting of the profits and weights of the items, and \textit{group-based} instance features. The latter include the profit-to-weight ratio, sum, minimum and maximum of the profits and weights for each group. Furthermore, for each group, we consider the relative share of its profit-to-weight ratio, total weight, and total profit relative to the sum across all groups. The generated data set is split into $N$ training and 100 test scenarios.
Every data point shown in the following figures is the average of 25 runs using different instances.

The implementation was carried out in Python 3.12.3 using Gurobi
version 13.0.2 as solver.
We use the CART and Lasso implementations of scikit-learn version 1.9.0.
Suitable values for the regularization parameter were determined empirically.
Up to 50 instances were solved in parallel with Gurobi's time limit set to one hour and its thread parameter to one. The code and data are available on GitHub\footnote{\url{https://github.com/sbstnmrtn/linear_interpretable_surrogates}}.

\begin{figure}[htbp]
    \centering
    \begin{subfigure}[t]{0.49\textwidth}
         \centering
         \includegraphics[width=\linewidth]{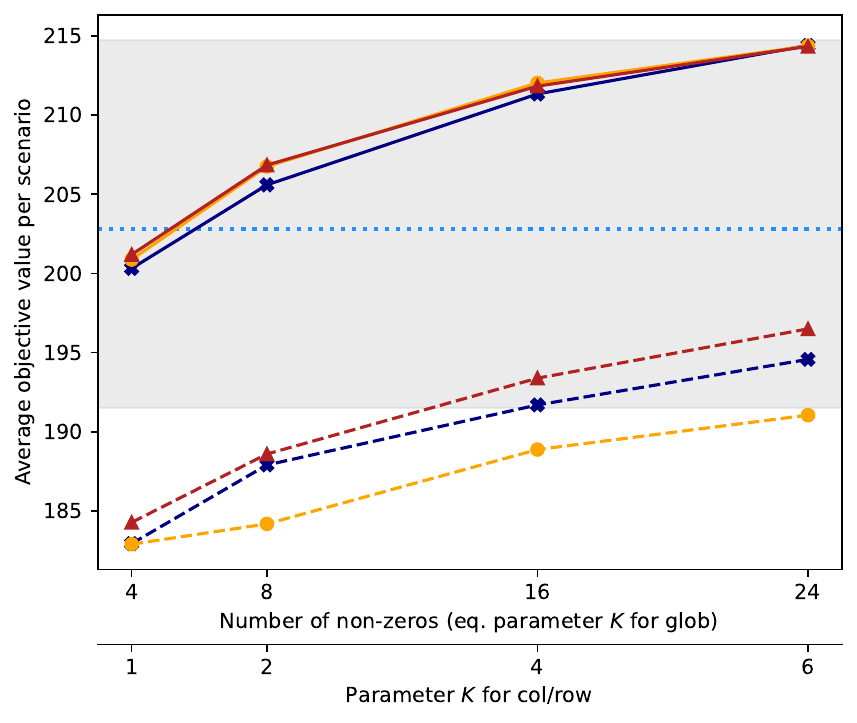}
         \caption{Evaluating on training data.}
         \label{fig:obj_k:tr}
    \end{subfigure}
    \begin{subfigure}[t]{0.49\textwidth}
         \centering
         \includegraphics[width=\linewidth]{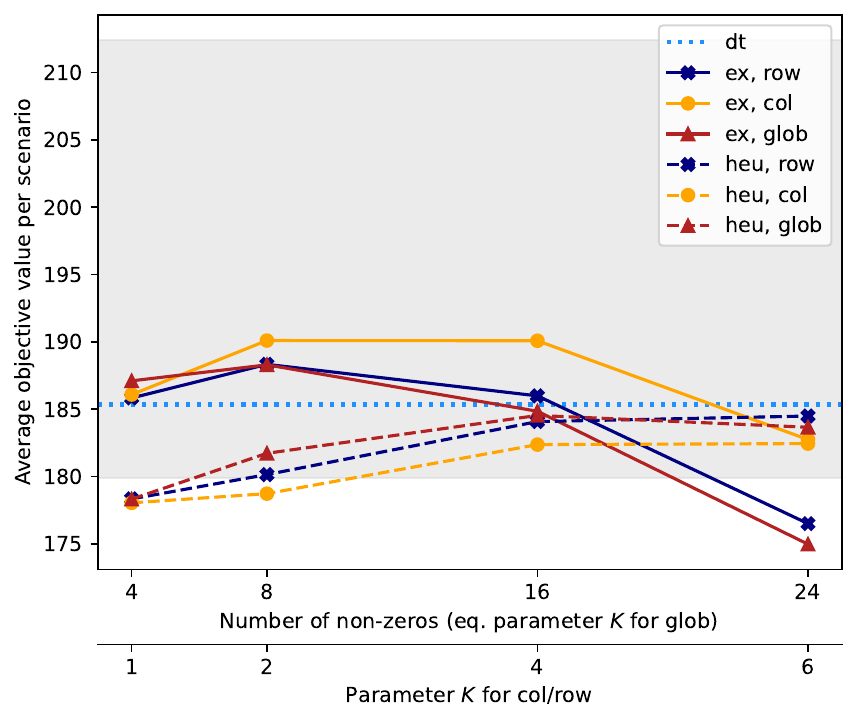}
         \caption{Evaluating on test data.}
         \label{fig:obj_k:te}
    \end{subfigure}
    \caption{Objective value vs. number of non-zeros using $N:= 10$. The legend in the right-hand plot applies to both figures.}
    \label{fig:obj_k}
\end{figure}

In Figure~\ref{fig:obj_k} the relationship between different values of the sparsity parameter $K$ and the objective value is shown. The gray area indicates the range given by the objective which would be achieved by solving every instance to optimality (upper bound) and using a single meta-solution (lower bound).

It can be seen that, as expected, the use of less sparse---and therefore less interpretable---models increases the performance on the training data. For $K=24$, the exact method even nearly yields optimal solutions. Nevertheless, the exact method shows clear indications of overfitting, while column sparsity makes these models more robust to overfitting. In most cases the exact method performs better than the heuristic, although the latter seems to be more resistant to overfitting. 
While the exact method reaches the time limit in most cases, the decision trees and the heuristic regression models were generated in a fraction of a second each.

\begin{figure}[htbp]
    \centering
    \begin{subfigure}[t]{0.49\textwidth}
         \centering
         \includegraphics[width=\linewidth]{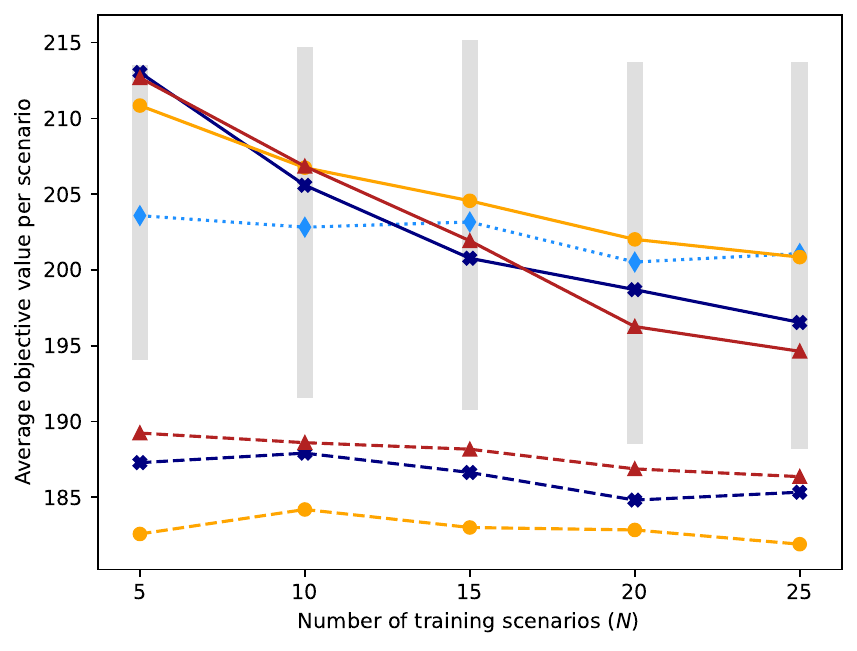}
         \caption{Evaluating on training data.}
         \label{fig:obj_n:tr}
    \end{subfigure}
    \begin{subfigure}[t]{0.49\textwidth}
         \centering
         \includegraphics[width=\linewidth]{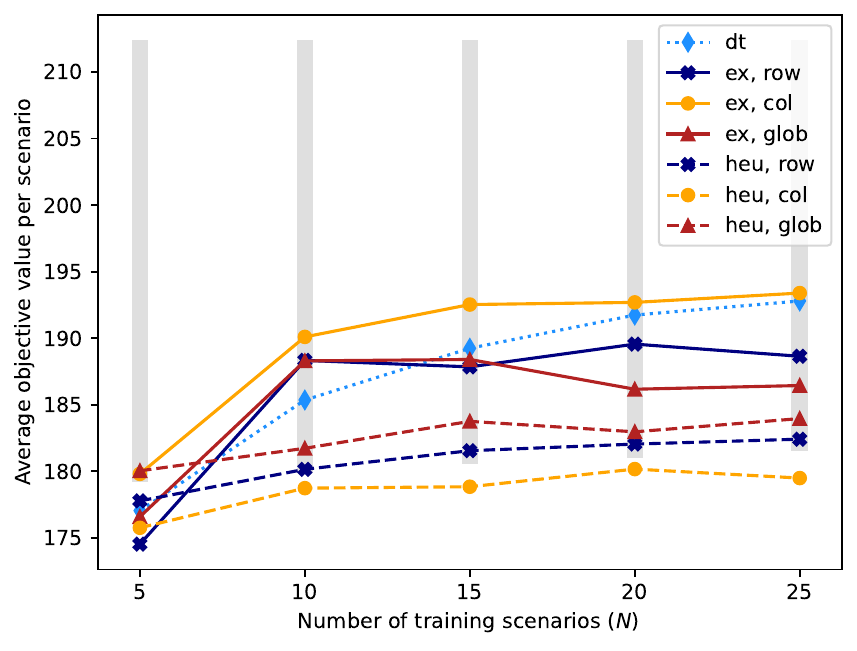}
         \caption{Evaluating on test data.}
         \label{fig:obj_n:te}
    \end{subfigure}
    \caption{Objective value vs. number of training scenarios using $K:=8$ for glob and $K:=2$ for col/row. The legend in the right-hand plot applies to both figures.}
    \label{fig:obj_n}
\end{figure}

In Figure~\ref{fig:obj_n} experiments varying the number of training scenarios are shown. Increasing $N$ generally leads to a decrease in the objective value per scenario on the training data, but to an increase on the unseen test data. In particular, an improvement in the quality of the decisions trees and of the exact method using column-based sparsity can be observed. The heuristic performs noticeably worse in the latter setting, which can be attributed to a misalignment of the column-based sparsity definition and Lasso.
Overall, sparse linear surrogates can outperform decision trees and may be particularly useful when linear explanations are more accessible to the target audience.

\section{Conclusion and Outlook}
\label{sec:outlook}

In this work, we proposed an interpretable optimization approach for budget-allocation problems using novel sparse linear surrogates. To compute them, we presented a mathematical programming formulation, which provides potentially exact solutions but is computationally demanding, and a heuristic based on out-of-the-box algorithms offering a favorable trade-off between runtime and solution quality. Our computational experiments demonstrate that these surrogates are a viable alternative to the established decision-tree-based surrogates. Future research may include the development of further algorithms to compute these surrogates as well as their application to other problem types and domains.

\end{document}